\documentclass[11pt]{amsart}
\usepackage{amsmath,amssymb,amsthm,esint}
\usepackage[a4paper,margin=2.7cm]{geometry}
\usepackage[hidelinks]{hyperref}
\hypersetup{pdftitle={A Garnett-Jones estimate for matrix weights}}
\newcommand{\R}{\mathbb R}
\newcommand{\T}{\mathbb T}
\newcommand{\B}{\mathcal B}
\DeclareMathOperator{\dist}{dist}
\newtheorem{theorem}{Theorem}[section]
\newtheorem{lemma}[theorem]{Lemma}
\newtheorem{proposition}[theorem]{Proposition}
\newtheorem{corollary}[theorem]{Corollary}
\newtheorem{introthm}{Theorem}

\theoremstyle{definition}
\newtheorem{definition}[theorem]{Definition}
\theoremstyle{remark}
\newtheorem{remark}[theorem]{Remark}
\title[A Garnett--Jones estimate for matrix weights]{A Garnett--Jones estimate for matrix weights}
\begin{document}
\author[M.\ Nielsen]{Morten Nielsen}
\author[H.\ \v{S}iki\'c]{Hrvoje \v{S}iki\'c}
%\thanks{The second named author has been partially
% supported by \ldots}
\address{
Department of Mathematical Sciences\\ Aalborg University\\ Thomas Manns Vej 23\\DK~-~9220 Aalborg\\ Denmark
}
\email{mnielsen@math.aau.dk}
\address{Department of Mathematics\\ Faculty of Science\\ University of Zagreb, Bijeni\v{c}ka 30\\ HR-10000\\
Zagreb, Croatia}
\email{hsikic@math.hr}
\begin{abstract}
For an \(N\times N\) matrix \(A_p\) weight \(W\) on an averaging
basis, set \(M_p=\max\{1,p-1\}\), \(1\leq p<\infty\).
At rate \(1/M_p\), the operator-norm exponential moment of
\(\log W\) centered on each basis set is bounded by
\(2[W]_{\mathbf A_p}^{1/M_p}\).
With scalar Garnett--Jones constant \(G\), the BMO distance
from \(\log W\) to bounded Hermitian functions is at most
\((2N^2-N)GM_p\), uniformly in the characteristic.
Power weights on intervals show the sharp order in \(p\).
Conversely, a uniform distance bound for scalar \(A_2\) logarithms
gives the scalar property with the same constant. Under a uniform
small-BMO condition this is equivalent to scalar \(A_2\) decomposition.
The estimates apply to cubes and axis-parallel rectangles on
\(\R^d\) and \(\T^d\). For rectangles with \(d>1\), the scalar
quadratic estimate bounds the distance by
\(\sqrt{2N^2-N}\,C_d M_p\), where \(C_d\) depends only on \(d\).
The rectangular BMO distance is comparable to the critical
exponential scale.
\end{abstract}
\maketitle

\section{Introduction}\label{sec:introduction}
The Garnett--Jones theorem \cite{GarJon78a} (see also Varopoulos
\cite{Var80a}) relates the BMO distance of a scalar function to
\(L^\infty\) to its exponential mean oscillation. If
\(\sup_Q\fint_Q e^{a|f-f_Q|}<\infty\) over all cubes
\(Q\subset\R^d\), then \(\dist_{\mathrm{BMO}}(f,L^\infty)\leq G_d/a\),
where \(G_d\) depends only on \(d\) and not on the size of the
exponential integral. The link with weights is that
\(w\in A_2\) forces exponential oscillation of \(\log w\) at
rate one, so \(\dist_{\mathrm{BMO}}(\log w,L^\infty)\leq G_d\)
uniformly over \(A_2\).

Nielsen and \v{S}iki\'c proved the scalar theorem for all
axis-parallel rectangles in each fixed dimension \(d>1\)
\cite[Theorem 1.1]{NieSik26}. The scalar space is little bmo,
with mean oscillation computed about the average on each rectangle.
The proof's quadratic oscillation bound gives the matrix distance
estimate for rectangles below.

This paper extends this picture to matrix weights. For a matrix
weight \(W\), the logarithm \(\log W\) is a Hermitian matrix
function, and we measure its mean oscillation in the Euclidean
operator norm. The basis \(\B\) is a general averaging basis on
\(\Omega=\R^d\) or \(\Omega=\T^d\). We use Roudenko's
mixed-integral characteristic for \(p>1\) and the
Frazier--Roudenko characteristic for \(p=1\)
(Section~\ref{sec:background}).
Throughout, \(M_p=\max\{1,p-1\}\).

\begin{introthm}\label{thmA}
Let \(1\leq p<\infty\) and \(W\in\mathbf A_p(\B)\). Then
\[
 \sup_{E\in\B}\fint_E
 \exp\left(\frac{\|\log W-(\log W)_E\|}{M_p}\right)
 \leq 2[W]_{\mathbf A_p(\B)}^{1/M_p}.
\]
\end{introthm}

The prefactor \(2\) is independent of the matrix size, the
dimension, and the geometry of \(\B\). The proof rests on a two-point comparison
of \(W(x)\) and \(W(y)\) in the Loewner order together with operator
monotonicity of the logarithm; no commutation between \(W(x)\) and
\(W(y)\) is needed.

We say that \(\B\) has the \emph{scalar Garnett--Jones property}
with constant \(G>0\) if, for every real \(f\in\mathrm{BMO}(\B)\)
and every \(a>0\),
\begin{equation}\label{eq:scalar-gj-basis}
 \sup_{E\in\B}\fint_E e^{a|f-f_E|}<\infty
 \quad\Longrightarrow\quad
 \dist_{\mathrm{BMO}(\B)}(f,L^\infty)\leq\frac{G}{a}.
\end{equation}
The optimal constant \(G^{\mathrm{opt}}_\B\) is the infimum of all
such \(G\), with \(G^{\mathrm{opt}}_\B=+\infty\) if there is none.
This property is an additional assumption on \(\B\); only the
upper distance estimate in \eqref{eq:scalar-gj-basis} is required.
For a weight \(W\) with locally integrable logarithm, write
\[
 \mathcal D_\B(W)=\dist_{\mathrm{BMO}_{\mathrm{op}}(\B)}
 \bigl(\log W,L^\infty(\Omega;\mathcal H_N)\bigr),
\]
where \(\mathcal H_N\) is the space of Hermitian \(N\times N\)
matrices.

\begin{introthm}\label{thmB}
If \(\B\) has the scalar Garnett--Jones property with constant
\(G\), then every \(N\times N\) weight \(W\in\mathbf A_p(\B)\),
\(1\leq p<\infty\), satisfies
\[
 \mathcal D_\B(W)\leq (2N^2-N)\,G\,M_p .
\]
For cubes and axis-parallel rectangles on \(\R^d\) and \(\T^d\),
the scalar constant \(G\) can be chosen to depend only on \(d\).
For the rectangle basis \(\mathcal R\) in dimension
\(d>1\), the upper bound improves to
\[
 \mathcal D_{\mathcal R}(W)
 \leq \sqrt{2N^2-N}\,C_d M_p,
\]
where \(C_d\) is the real-valued correction constant of
\cite[Theorem 1.1 and Section 8.1]{NieSik26}.
\end{introthm}

The infimum defining \(\mathcal D_\B(W)\) can be taken over
approximants \(\log V\) with both \(V\) and \(V^{-1}\) bounded.
The distance bound is independent of the \(A_p\) characteristic
because the exponential rate in Theorem~\ref{thmA} is uniform.
The factor \(2N^2-N\) comes from summing the scalar coordinate
errors. On rectangles the quadratic estimate lets us sum their
squares. The order in \(p\) is sharp already on intervals, where
scalar power weights give a lower bound of \(\tfrac12 M_p\log2\)
for the supremum of the distance over matrix \(A_p\) weights of
each fixed size
(Corollary~\ref{cor:sharp-p}).

Scalar \(A_2\) weights detect the necessity of the scalar
hypothesis in Theorem~\ref{thmB}. Set
\[
 L_{\B,N}=\sup\{\mathcal D_\B(W):
 W\in\mathbf A_2(\B),\ W\text{ is }N\times N\}.
\]

\begin{introthm}\label{thmC}
For every averaging basis \(\B\) and every \(N\geq1\), as extended
real numbers,
\[
 G^{\mathrm{opt}}_\B=L_{\B,1}\leq L_{\B,N}\leq(2N^2-N)G^{\mathrm{opt}}_\B .
\]
If moreover \(\B\) satisfies \eqref{eq:small-bmo-exponential},
which uniformly controls the \(A_2\) characteristic of \(e^h\)
for real \(h\) of small BMO seminorm, then the scalar Garnett--Jones
property, finiteness of \(L_{\B,N}\), and the \(A_2\)-decomposition
property of \cite{NieSik19} are equivalent.
\end{introthm}

\subsection*{Relation to earlier work}
Bloom proved that \(\log W\in\mathrm{BMO}\) on the circle if the
Hilbert transform is bounded on \(L^2(W)\)
\cite[Theorem 5.4]{Blo81}. This boundedness characterizes matrix
\(A_2\) by Treil and Volberg \cite{TreVol97}. Bloom also constructed
a Hermitian BMO function whose exponential falls outside matrix
\(A_2\) at every positive rate \cite[pp.~100--103]{Blo81}.

The two-sided estimates in \cite{NieSik21,NieSik25} use a
prescribed log-preserving unitary matrix function \(U(x)\) for
weights \(U^*e^DU\). The distance in \cite[Definition 17.6]{NieSik25},
for example, measures \(D-D'\) in diagonal coordinates and uses
approximants of the form \(U^*D'U\) with \(D'\) bounded and
diagonal. Here the distance uses the operator norm and allows
every bounded Hermitian approximant. Section~\ref{sec:comparison}
relates this distance to the matrix \(A_2\) threshold.

In \cite{NieSik19}, the scalar \(A_2\)-decomposition property was
characterized on structured Muckenhoupt bases with the
Calder\'on--Zygmund decomposition property by an upper distance bound
in terms of the scalar \(A_2\) threshold
\(\varepsilon_\B(f)=\inf\{\lambda>0:e^{f/\lambda}\in A_2(\B)\}\).
We identify that bound with the scalar Garnett--Jones property,
with the same constant, and isolate the small-BMO condition for
the converse from distance to decomposition. The scalar rectangle
theorem and its John--Nirenberg estimate give the scalar
\(A_2\)-decomposition property for rectangles
(Corollary~\ref{cor:rect-decomposition}).

\subsection*{Organization}
Section~\ref{sec:background} fixes definitions and introduces the
scalar John--Nirenberg property of a basis.
Section~\ref{sec:logarithms} proves Theorem~\ref{thmA}.
Section~\ref{sec:distance} proves Theorem~\ref{thmB} for general
bases, gives the cube and rectangle applications on Euclidean
space and the torus, and proves a two-sided estimate under a
John--Nirenberg hypothesis. Section~\ref{sec:equivalence} proves
Theorem~\ref{thmC}. Section~\ref{sec:sharpness} shows that the
factor \(M_p\) is sharp in order. Section~\ref{sec:comparison}
compares with \cite{NieSik21,NieSik25}, and
Section~\ref{sec:questions} lists open questions. Apart from the
order of growth in \(p\), the constants are not claimed to be
optimal.

\section{Averaging bases and matrix weights}\label{sec:background}
Let \(\Omega\in\{\R^d,\T^d\}\), where
\(\T^d=\R^d/\mathbb Z^d\), with Lebesgue measure. An
\emph{averaging basis} \(\B\) is a family of relatively compact
measurable sets \(E\subset\Omega\) with \(0<|E|<\infty\), such
that every compact subset of \(\Omega\) is covered, up to a null
set, by finitely many members of \(\B\). We write
\(f_E=\fint_E f=|E|^{-1}\int_E f\).
No differentiation, doubling, or covering-lemma
assumption is imposed.

The \emph{cube basis} \(\mathcal Q\) on \(\R^d\) consists of all
axis-parallel cubes. On \(\T^d\), it consists of the images of
axis-parallel cubes of side length at most one, including those
crossing the boundary of the fundamental domain; the cubes of side
one give the whole torus.

The \emph{rectangle basis} \(\mathcal R\) on \(\R^d\) consists
of all products of bounded intervals of positive length. On
\(\T^d\) it consists of the images of those products whose side
lengths are at most one, including rectangles that cross the
boundary of the fundamental domain. A side of length one gives a
full circle in that coordinate. When \(d=1\), the cube and
rectangle bases coincide.

A \emph{matrix weight} is a measurable, locally integrable
\(N\times N\) matrix function that is Hermitian and positive
definite almost everywhere. We write \(\|\cdot\|\) for the
Euclidean operator norm and \(\mathcal H_N\) for the real vector space of
Hermitian \(N\times N\) matrices. For \(1<p<\infty\), put
\(p'=p/(p-1)\) and
\begin{equation}\label{eq:ap-definition}
 [W]_{\mathbf A_p(\B)}
 =\sup_{E\in\B}\fint_E
 \left(\fint_E\|W^{1/p}(x)W^{-1/p}(y)\|^{p'}dy\right)^{p/p'}dx,
\end{equation}
the formulation of matrix \(A_p\) due to Roudenko
\cite[Lemma 1.3]{Rou03a}, used here as the definition on \(\B\).
At the endpoint we use the matrix \(A_1\) condition of Frazier and
Roudenko \cite[(1.2) with \(p=1\)]{FraRou04},
\begin{equation}\label{eq:a1-definition}
 [W]_{\mathbf A_1(\B)}
 =\sup_{E\in\B}\operatorname*{ess\,sup}_{y\in E}
      \fint_E\|W(x)W^{-1}(y)\|\,dx.
\end{equation}
The class \(\mathbf A_p(\B)\) consists of matrix weights with
finite characteristic. For \(N=1\) these reduce to the usual
scalar characteristics
\[
 [w]_{A_p(\B)}
 =\sup_E\Bigl(\fint_Ew\Bigr)
          \Bigl(\fint_Ew^{-1/(p-1)}\Bigr)^{p-1}\ (p>1),
 \qquad
 [w]_{A_1(\B)}
 =\sup_E\operatorname*{ess\,sup}_{y\in E}
           \frac{\fint_Ew}{w(y)}.
\]

For \(H\in L^1_{\mathrm{loc}}(\Omega;\mathcal H_N)\) and scalar
\(f\in L^1_{\mathrm{loc}}(\Omega)\), set
\begin{equation}\label{eq:matrix-bmo}
 \|H\|_{\mathrm{BMO}_{\mathrm{op}}(\B)}
 =\sup_{E\in\B}\fint_E\|H-H_E\|,
 \qquad
 \|f\|_{\mathrm{BMO}(\B)}
 =\sup_{E\in\B}\fint_E|f-f_E|.
\end{equation}
Since taking an entry commutes with averaging and
\(|X_{ij}|\leq\|X\|\leq\sum_{i,j}|X_{ij}|\),
\begin{equation}\label{eq:entrywise}
 \max_{i,j}\|H_{ij}\|_{\mathrm{BMO}(\B)}
 \leq\|H\|_{\mathrm{BMO}_{\mathrm{op}}(\B)}
 \leq\sum_{i,j}\|H_{ij}\|_{\mathrm{BMO}(\B)}.
\end{equation}

Distances are infima of the seminorm of the difference:
\[
 \dist_{\mathrm{BMO}(\B)}(f,L^\infty)
 =\inf_{b\in L^\infty(\Omega;\R)}\|f-b\|_{\mathrm{BMO}(\B)}
\]
for real \(f\). The matrix distance is defined in the same way
using bounded Hermitian approximants. The seminorms and
distances may take the value \(+\infty\) for locally integrable
functions. For a matrix weight \(W\) with locally integrable
logarithm, write
\begin{equation}\label{eq:distance-definition}
 \mathcal D_{\B}(W)=\dist_{\mathrm{BMO}_{\mathrm{op}}(\B)}
       \bigl(\log W,L^\infty(\Omega;\mathcal H_N)\bigr).
\end{equation}
The logarithm is defined by functional calculus. Lemma~\ref{lem:integrability}
gives its local integrability for every \(W\in\mathbf A_p(\B)\),
\(1\leq p<\infty\).

\begin{remark}\label{rem:weight-approximants}
Let \(\mathcal V_N\) be the class of matrix weights \(V\) with
\(c_VI_N\leq V\leq C_VI_N\) a.e.\ for some
\(0<c_V\leq C_V<\infty\). Then
\(\mathcal D_\B(W)=\inf_{V\in\mathcal V_N}
\|\log W-\log V\|_{\mathrm{BMO}_{\mathrm{op}}(\B)}\).
Indeed, \(V\mapsto\log V\) maps \(\mathcal V_N\) onto
\(L^\infty(\Omega;\mathcal H_N)\): if \(V\in\mathcal V_N\) then
\(\|\log V\|_\infty\leq\max\{|\log c_V|,|\log C_V|\}\), and if
\(B\) is bounded Hermitian then
\(e^{-\|B\|_\infty}I_N\leq e^B\leq e^{\|B\|_\infty}I_N\).
\end{remark}

In addition to the scalar Garnett--Jones property
\eqref{eq:scalar-gj-basis}, we use the following scalar property.

\begin{definition}\label{def:scalar-jn}
The basis \(\B\) has the \emph{scalar John--Nirenberg property}
with constants \(J\geq1\) and \(c>0\) if every real
\(h\in\mathrm{BMO}(\B)\) of positive seminorm satisfies
\begin{equation}\label{eq:scalar-john-nirenberg}
 \frac{|\{x\in E:|h(x)-h_E|>t\}|}{|E|}
 \leq J\exp\left(-\frac{ct}{\|h\|_{\mathrm{BMO}(\B)}}\right)
 \qquad(E\in\B,\ t>0).
\end{equation}
\end{definition}

The John--Nirenberg tail estimate gives the following uniform
exponential moment bound.

\begin{lemma}\label{lem:layer-cake}
Let \(\B\) have the scalar John--Nirenberg property with
constants \(J,c\). Let \(s>0\) and let \(h\in\mathrm{BMO}(\B)\) be
real with \(\|h\|_{\mathrm{BMO}(\B)}\leq s\). Then, for every
\(E\in\B\) and \(0<a<c/s\),
\[
 \fint_E e^{a|h-h_E|}\leq1+\frac{Ja}{c/s-a}.
\]
\end{lemma}

\begin{proof}
If \(\|h\|_{\mathrm{BMO}(\B)}=0\), then \(h=h_E\) a.e.\ on \(E\)
and the average is one. Otherwise
\eqref{eq:scalar-john-nirenberg} gives the tail bound
\(J e^{-ct/s}\), and the layer-cake formula yields
\[
 \fint_E e^{a|h-h_E|}
 =1+\int_0^\infty ae^{at}\,
       \frac{|\{x\in E:|h(x)-h_E|>t\}|}{|E|}\,dt
 \leq1+Ja\int_0^\infty e^{-(c/s-a)t}\,dt.
 \qedhere
\]
\end{proof}

\section{Exponential oscillation of the logarithm}\label{sec:logarithms}

To justify averaging \(\log W\), we first control the relevant
inverse powers of \(W\). For \(p>1\), local integrability of
\(W^{-1/(p-1)}\) follows from Roudenko's matrix \(A_p\) duality
for Euclidean balls \cite[Corollary 3.3]{Rou03a}. For structured
torus bases see \cite[Definition 17.2 and Remark 17.1]{NieSik25}.
At \(p=1\), local essential boundedness of \(W^{-1}\) follows
from the reducing-operator estimate on dyadic cubes
in \cite[Lemma 5.4]{FraRou04}. Cruz-Uribe, Moen, and Rodney give
both bounds for Euclidean domains
\cite[Lemma 4.5 and proof of Theorem 5.2]{CruMoeRod16}.
We prove the bounds directly for general averaging bases and
deduce local integrability of \(\log W\).

\begin{lemma}\label{lem:integrability}
If \(W\in\mathbf A_p(\B)\), \(1<p<\infty\), then
\(W^{-1/(p-1)}\) is locally integrable. If
\(W\in\mathbf A_1(\B)\), then \(W^{-1}\) is locally
essentially bounded. In either case
\(\log W\in L^1_{\mathrm{loc}}(\Omega;\mathcal H_N)\).
\end{lemma}

\begin{proof}
Fix \(E\in\B\). For \(p>1\), finiteness of
\eqref{eq:ap-definition} lets us choose \(x\in E\) with \(W(x)\) positive definite
and \(y\mapsto\|W^{1/p}(x)W^{-1/p}(y)\|\) in \(L^{p'}(E)\). From
\[
 \|W^{-1/p}(y)\|
 \leq\|W^{-1/p}(x)\|\,\|W^{1/p}(x)W^{-1/p}(y)\|
\]
we get \(W^{-1/p}\in L^{p'}(E)\), and
\(\|W^{-1/p}\|^{p'}=\|W^{-1/(p-1)}\|\) gives integrability of the
dual power on \(E\). Finite coverings of compact sets give local
integrability.

Let \(p=1\). For almost every \(x,y\in E\),
\(\|W^{-1}(y)\|\leq\|W^{-1}(x)\|\,\|W(x)W^{-1}(y)\|\), that is,
\[
 \|W^{-1}(y)\|\,\|W^{-1}(x)\|^{-1}\leq\|W(x)W^{-1}(y)\|.
\]
Averaging in \(x\) gives, for almost every \(y\in E\),
\[
 \|W^{-1}(y)\|
 \fint_E\|W^{-1}(x)\|^{-1}dx
 \leq[W]_{\mathbf A_1(\B)}.
\]
The integrand \(\|W^{-1}(x)\|^{-1}=\lambda_{\min}(W(x))\) is
positive a.e.\ and bounded by \(\|W(x)\|\), so the average is
positive and finite, and \(W^{-1}\) is essentially bounded on
\(E\). Again finite coverings conclude.

For the logarithm, \(\log t\leq t\) for \(t\geq1\) and
\(-\log t\leq(p-1)t^{-1/(p-1)}\) for \(0<t<1\) give
\(\|\log W\|\leq\|W\|+(p-1)\|W^{-1/(p-1)}\|\) when \(p>1\).
When \(p=1\), \(\|W\|\) controls the positive eigenvalues of
\(\log W\) and the local bound on \(W^{-1}\) the negative ones.
\end{proof}

To control the oscillation of \(\log W\), we relate its two-point
differences to the matrix products appearing in the \(A_p\)
characteristic. For \(p\geq1\) and \(x,y\) with \(W(x),W(y)\)
positive definite, set
\[
 C_p(x,y)=\|W^{1/p}(x)W^{-1/p}(y)\|.
\]
The two constants \(C_p(x,y)\) and \(C_p(y,x)\) are the norms of
a matrix and of its inverse; they control the upper and the lower
part of the spectrum of \(\log W(x)-\log W(y)\), respectively.

\begin{lemma}\label{lem:two-point}
Let \(W\) be a matrix weight, \(H=\log W\), and \(p\geq1\). For
almost every \(x,y\),
\begin{equation}\label{eq:log-pair}
 \|H(x)-H(y)\|
 \leq p\log\max\{C_p(x,y),C_p(y,x)\}.
\end{equation}
\end{lemma}

\begin{proof}
Fix \(x,y\) with \(W(x),W(y)\) positive definite and put
\(T=W^{1/p}(x)W^{-1/p}(y)\), so \(\|T\|=C_p(x,y)\). Since
\[
 T^*T=W^{-1/p}(y)W^{2/p}(x)W^{-1/p}(y)
 \leq C_p(x,y)^2I_N,
\]
congruence by \(W^{1/p}(y)\) gives
\[
 W^{2/p}(x)\leq C_p(x,y)^2\,W^{2/p}(y).
\]
The logarithm is operator monotone on positive definite matrices
\cite[Section 5.3.7]{Bha07}:
\[
 \log P=\int_0^\infty\bigl((1+t)^{-1}I_N-(P+tI_N)^{-1}\bigr)\,dt,
\]
and \(P\mapsto-(P+tI_N)^{-1}\) is operator monotone. Taking the
logarithm of the preceding Loewner inequality gives
\[
 \tfrac2pH(x)\leq 2\log C_p(x,y)\,I_N+\tfrac2pH(y).
\]
Interchanging \(x\) and \(y\),
\[
 -p\log C_p(y,x)\,I_N
 \leq H(x)-H(y)\leq p\log C_p(x,y)\,I_N.
\]
Now \(T^{-1}=W^{1/p}(y)W^{-1/p}(x)\), so
\(1\leq\|T\|\,\|T^{-1}\|=C_p(x,y)C_p(y,x)\). Hence
\[
 \theta=p\log\max\{C_p(x,y),C_p(y,x)\}\geq0
\]
dominates both \(p\log C_p(x,y)\) and \(p\log C_p(y,x)\), and
\(-\theta I_N\leq H(x)-H(y)\leq\theta I_N\). The eigenvalues of
the Hermitian matrix \(H(x)-H(y)\) therefore lie in
\([-\theta,\theta]\), which is \eqref{eq:log-pair}.
\end{proof}

To use the two-point comparison, we convert the mixed \(A_p\)
characteristic into a double-integral moment bound with the same
exponent in both variables. Set \(m_p=p/M_p\), so that
\(m_p=\min\{p,p'\}\) for \(p>1\) and \(m_1=1\).

\begin{lemma}\label{lem:kernel-moment}
Let \(1<p<\infty\), \(W\in\mathbf A_p(\B)\), and
\(\kappa=[W]_{\mathbf A_p(\B)}\). Then, for every \(E\in\B\),
\begin{equation}\label{eq:kernel-moment}
 \fint_E\fint_E C_p(x,y)^{m_p}\,dy\,dx
 \leq \kappa^{1/M_p}.
\end{equation}
\end{lemma}

\begin{proof}
If \(p\leq2\), then \(m_p=p\leq p'\), \(M_p=1\), and Jensen's
inequality in \(y\) gives
\[
 \fint_E\fint_E C_p^p\,dy\,dx
 \leq\fint_E\Bigl(\fint_E C_p^{p'}\,dy\Bigr)^{p/p'}dx
 \leq\kappa .
\]
If \(p\geq2\), then \(m_p=p'\), \(M_p=p-1\), \(p'/p=1/(p-1)\), and
Jensen's inequality in \(x\) gives
\[
 \fint_E\fint_E C_p^{p'}\,dy\,dx
 \leq\Bigl(\fint_E\Bigl(\fint_E C_p^{p'}\,dy\Bigr)^{p/p'}dx\Bigr)^{p'/p}
 \leq\kappa^{1/(p-1)} .
\qedhere
\]
\end{proof}

We can now prove Theorem~\ref{thmA}, together with the resulting
BMO bound.

\begin{theorem}\label{thm:log-oscillation}
Let \(1\leq p<\infty\), \(W\in\mathbf A_p(\B)\), and
\(\kappa=[W]_{\mathbf A_p(\B)}\). Then
\(H=\log W\in L^q_{\mathrm{loc}}\) for every \(1\leq q<\infty\),
\begin{equation}\label{eq:log-exponential}
 \sup_{E\in\B}\fint_E
 \exp\left(\frac{\|H(x)-H_E\|}{M_p}\right)dx
 \leq 2\kappa^{1/M_p},
\end{equation}
and
\begin{equation}\label{eq:log-bmo}
 \|\log W\|_{\mathrm{BMO}_{\mathrm{op}}(\B)}
 \leq \log\kappa+M_p\log2.
\end{equation}
\end{theorem}

\begin{proof}
Fix \(E\in\B\). We first show
\begin{equation}\label{eq:log-double}
 \fint_E\fint_E
 \exp\left(\frac{\|H(x)-H(y)\|}{M_p}\right)dy\,dx
 \leq2\kappa^{1/M_p}.
\end{equation}
By Lemma~\ref{lem:two-point}, \(p/M_p=m_p\), and
\(\max\{s,t\}^{m_p}\leq s^{m_p}+t^{m_p}\),
\[
 \exp\left(\frac{\|H(x)-H(y)\|}{M_p}\right)
 \leq C_p(x,y)^{m_p}+C_p(y,x)^{m_p}.
\]
Averaging over \(E\times E\), the two terms have the same double
average. For \(p>1\) it is at most \(\kappa^{1/M_p}\) by
Lemma~\ref{lem:kernel-moment}. For \(p=1\) the definition
\eqref{eq:a1-definition} gives
\(\fint_E\fint_E C_1(x,y)\,dx\,dy\leq\kappa\) directly. This
proves \eqref{eq:log-double}.

We next pass from two-point differences to oscillation about the
mean. By Lemma~\ref{lem:integrability}, \(H\) is integrable on
\(E\), and for almost every \(x\in E\),
\[
 H(x)-H_E=\fint_E\bigl(H(x)-H(y)\bigr)\,dy.
\]
Taking norms and using the triangle inequality gives
\[
 \|H(x)-H_E\|\leq\fint_E\|H(x)-H(y)\|\,dy.
\]
Since the scalar exponential is increasing and convex, Jensen's
inequality in \(y\) now yields
\[
 \exp\left(\frac{\|H(x)-H_E\|}{M_p}\right)
 \leq\exp\left(\fint_E\frac{\|H(x)-H(y)\|}{M_p}\,dy\right)
 \leq\fint_E
 \exp\left(\frac{\|H(x)-H(y)\|}{M_p}\right)dy.
\]
Averaging in \(x\) and applying \eqref{eq:log-double} proves
\eqref{eq:log-exponential}.

To obtain the BMO bound, apply scalar Jensen once more, now to
\(x\mapsto\|H(x)-H_E\|/M_p\) on \(E\):
\[
 \exp\left(\frac1{M_p}\fint_E\|H(x)-H_E\|\,dx\right)
 \leq\fint_E\exp\left(\frac{\|H(x)-H_E\|}{M_p}\right)dx
 \leq2\kappa^{1/M_p}.
\]
Taking logarithms and multiplying by \(M_p\), we obtain
\[
 \fint_E\|H-H_E\|\leq M_p\log2+\log\kappa.
\]
The supremum over \(E\in\B\) is \eqref{eq:log-bmo}.
Finally, the exponential bound implies \(H-H_E\in L^q(E)\) for
every finite \(q\geq1\). Since \(H_E\) is a fixed finite matrix,
\(H\in L^q(E)\), and finite coverings of compact sets give the
local statement.
\end{proof}

\section{The matrix Garnett--Jones estimate}\label{sec:distance}

We obtain a bounded Hermitian approximant by assembling scalar
approximants to the real coordinates of \(H\).

\begin{lemma}\label{lem:matrix-distance}
Let \(\B\) have the scalar Garnett--Jones property with constant
\(G\), and let \(H\in L^1_{\mathrm{loc}}(\Omega;\mathcal H_N)\).
If, for some \(a>0\),
\[
 \sup_{E\in\B}\fint_E e^{a\|H-H_E\|}<\infty,
\]
then
\begin{equation}\label{eq:matrix-distance}
 \dist_{\mathrm{BMO}_{\mathrm{op}}(\B)}
 \bigl(H,L^\infty(\Omega;\mathcal H_N)\bigr)
 \leq\frac{(2N^2-N)G}a.
\end{equation}
\end{lemma}

\begin{proof}
Consider the \(N^2\) real functions \(H_{ii}\),
\(\operatorname{Re}H_{ij}\), and \(\operatorname{Im}H_{ij}\),
\(i<j\). The oscillation of each on \(E\) is at most
\(\|H-H_E\|\), so each has an exponential-oscillation bound at rate
\(a\) and, by Jensen, lies in \(\mathrm{BMO}(\B)\). Given
\(\varepsilon>0\), \eqref{eq:scalar-gj-basis} gives a bounded real
approximant to each with error at most \(G/a+\varepsilon\) in
\(\mathrm{BMO}(\B)\). Form a bounded Hermitian \(B\) using the
approximants to \(H_{ii}\) on the diagonal. For \(i<j\), let
\(B_{ij}=b^{\mathrm{re}}_{ij}+i b^{\mathrm{im}}_{ij}\), with
the two approximants to the real and imaginary parts of \(H_{ij}\),
and set \(B_{ji}=\overline{B_{ij}}\). Put \(R=H-B\).
By \eqref{eq:entrywise},
\(\|R\|_{\mathrm{BMO}_{\mathrm{op}}}\) is at most the sum of the
entrywise seminorms. Each diagonal entry contributes one error.
Each of the \(\binom N2\) off-diagonal pairs contributes four:
real and imaginary parts, in both positions \((i,j)\) and \((j,i)\).
Hence
\[
 \|R\|_{\mathrm{BMO}_{\mathrm{op}}(\B)}
 \leq \Bigl(N+4\tbinom N2\Bigr)\Bigl(\frac Ga+\varepsilon\Bigr)
 =(2N^2-N)\Bigl(\frac Ga+\varepsilon\Bigr).
\]
Let \(\varepsilon\downarrow0\).
\end{proof}

Applying the lemma to \(\log W\) gives the general-basis
estimate in Theorem~\ref{thmB}.

\begin{theorem}\label{thm:log-distance}
Let \(\B\) have the scalar Garnett--Jones property with constant
\(G\). If \(W\in\mathbf A_p(\B)\), \(1\leq p<\infty\), then
\begin{equation}\label{eq:log-distance}
 \mathcal D_{\B}(W)\leq (2N^2-N)\,G\,M_p.
\end{equation}
By Remark~\ref{rem:weight-approximants}, the bounded approximants
may be taken of the form \(\log V\) with \(V\in\mathcal V_N\).
\end{theorem}

\begin{proof}
By Theorem~\ref{thm:log-oscillation}, \(H=\log W\) satisfies the
hypothesis of Lemma~\ref{lem:matrix-distance} with \(a=1/M_p\).
\end{proof}

\subsection{Transfer to the torus}\label{subsec:periodization}

Write \(\mathcal A(\Omega)\) for either the cube or the rectangle
basis on \(\Omega\), and set
\[
 K_{a,\mathcal A(\Omega)}(f)
   =\sup_{E\in\mathcal A(\Omega)}\fint_E e^{a|f-f_E|}.
\]
For real \(f\in L^1(\T^d)\), we first control the exponential
moments of its periodic extension \(\widetilde f\) on \(\R^d\).
We then periodize a bounded Euclidean correction, preserving
its bounds for mean absolute or quadratic oscillation.

Fix a Euclidean test rectangle \(R=\prod_{i=1}^d I_i\), with
side lengths \(\ell_i=|I_i|\). Let \(S\subset\T^d\) be the
product of the images of the \(I_i\) with \(\ell_i\leq1\) and
the full circles in the other coordinates. For a cube \(R\),
\(S\) belongs to the torus cube basis. For any nonnegative
measurable \(g\) on \(S\) and its periodic lift to \(R\),
\begin{equation}\label{eq:periodic-domination}
 \fint_R \widetilde g\leq 2^d\fint_S g.
\end{equation}
Indeed, each interval of length \(\ell_i>1\) is covered by
\(\lceil\ell_i\rceil\) unit intervals, and
\(\lceil\ell_i\rceil/\ell_i\leq2\). The shorter sides project
injectively up to null sets. Apply \eqref{eq:periodic-domination}
to \(g=e^{a|f-f_S|}\). Jensen bounds
\(e^{a|\widetilde f_R-f_S|}\) by the same average on \(R\), so
\begin{equation}\label{eq:periodic-moment}
 K_{a,\mathcal A(\R^d)}(\widetilde f)
 \leq 4^d K_{a,\mathcal A(\T^d)}(f)^2.
\end{equation}
Applying \eqref{eq:periodic-domination} to \(|f-f_S|\) and
centering on \(R\) also gives
\begin{equation}\label{eq:periodic-bmo}
 \|\widetilde f\|_{\mathrm{BMO}(\mathcal A(\R^d))}
 \leq 2^{d+1}\|f\|_{\mathrm{BMO}(\mathcal A(\T^d))}.
\end{equation}

\begin{lemma}\label{lem:periodic-gj}
Let \(\mathcal A\) be the cube or rectangle basis and \(q\in\{1,2\}\).
Suppose the periodic extension \(\widetilde f\) of a real torus
function has a bounded real approximant \(b\) on \(\R^d\) with
\[
 \sup_{E\in\mathcal A(\R^d)}
 \left(\fint_E|(\widetilde f-b)-(\widetilde f-b)_E|^q\right)^{1/q}
 \leq\sigma.
\]
Then there is a real approximant on \(\T^d\) with the same bound
\(\sigma\) on the torus basis and \(L^\infty\) norm at most
\(\|b\|_\infty\). Consequently, the scalar Garnett--Jones property
passes from \(\mathcal A(\R^d)\) to \(\mathcal A(\T^d)\) with the
same constant.
\end{lemma}

\begin{proof}
For \(n\geq1\), average the integer translates of \(b\):
\[
 b_n(x)=\frac1{(2n+1)^d}\sum_{k\in\{-n,\ldots,n\}^d}b(x+k).
\]
Periodicity of \(\widetilde f\), translation invariance of the
Euclidean basis, and the triangle inequality in \(L^q(E)\) give
the same bound \(\sigma\) for \(\widetilde f-b_n\). Also
\(\|b_n\|_\infty\leq\|b\|_\infty\). In each coordinate direction
\(e_i\), the sum telescopes, giving
\[
 \|b_n(\cdot+e_i)-b_n\|_\infty
 \leq \frac{2\|b\|_\infty}{2n+1}.
\]
The uniform \(L^\infty\) bound and separability of \(L^1(\R^d)\)
give a weak-* convergent subsequence in \(L^\infty\).
Its limit \(b_{\mathrm{per}}\) is periodic because translation is
weak-* continuous. On each \(E\)
in the Euclidean basis the centered residuals converge weakly
in \(L^q(E)\). Indeed, \(\widetilde f\in L^q(E)\) by the
hypothesis, and the dual test functions, in \(L^\infty(E)\) for
\(q=1\) and \(L^2(E)\) for \(q=2\), extend by zero to \(L^1(\R^d)\).
Weak lower semicontinuity transfers the oscillation bound to
\(\widetilde f-b_{\mathrm{per}}\). Restricting to the Euclidean
lifts of the torus sets proves the assertion about the
approximant.

For the last assertion use \eqref{eq:periodic-moment} at the
given rate \(a\), apply the Euclidean scalar theorem with error
\(G/a+\delta\), and let \(\delta\downarrow0\).
\end{proof}

\subsection{Cubes and rectangles}

The scalar Garnett--Jones theorem \cite{GarJon78a}, see also
\cite[Introduction, Theorem (G.J.)]{Var80a}, gives the property
\eqref{eq:scalar-gj-basis} on Euclidean cubes with a constant
\(G_d\) depending only on \(d\).

\begin{corollary}\label{cor:cube-log-distance}
Let \(\mathcal Q\) be the cube basis on \(\R^d\) or \(\T^d\). If
\(W\in\mathbf A_p(\mathcal Q)\), \(1\leq p<\infty\), then
\(\mathcal D_{\mathcal Q}(W)\leq(2N^2-N)\,G_d\,M_p\).
\end{corollary}

\begin{proof}
Apply Theorem~\ref{thm:log-distance} with the classical scalar
theorem on \(\R^d\) and Lemma~\ref{lem:periodic-gj} on \(\T^d\).
\end{proof}

For Euclidean rectangles in dimension \(d>1\), the construction in
\cite[Section 8.1]{NieSik26} gives a real correction \(b\) for every
real \(f\in L^1_{\mathrm{loc}}(\R^d)\) with \(K_{a,\mathcal R}(f)<\infty\):
\begin{equation}\label{eq:scalar-rect-correction}
 \|b\|_\infty\leq a^{-1}\log K_{a,\mathcal R}(f),
 \qquad
 \sup_{R\in\mathcal R}
 \left(\fint_R|(f-b)-(f-b)_R|^2\right)^{1/2}
 \leq C_d/a .
\end{equation}
In particular, \(G=C_d\) works in \eqref{eq:scalar-gj-basis}
for Euclidean rectangles and, by Lemma~\ref{lem:periodic-gj},
for torus rectangles. The quadratic bound lets us sum the
squared errors of the real coordinates.

\begin{proposition}\label{prop:rect-correction}
Let \(d>1\) and let \(\mathcal R\) be the rectangle basis on
\(\Omega=\R^d\) or \(\T^d\). Suppose
\(H\in L^1_{\mathrm{loc}}(\Omega;\mathcal H_N)\) and
\[
 K=\sup_{R\in\mathcal R}\fint_R e^{a\|H-H_R\|}<\infty
 \qquad\text{for some }a>0 .
\]
Then there is a bounded Hermitian \(B\) such that
\begin{equation}\label{eq:rect-quadratic}
 \sup_{R\in\mathcal R}
 \left(\fint_R\|(H-B)-(H-B)_R\|^2\right)^{1/2}
 \leq \frac{\sqrt{2N^2-N}\,C_d}{a}.
\end{equation}
Moreover, \(\|B\|_\infty\leq a^{-1}\sqrt{2N^2-N}\log K^*\),
where \(K^*=K\) on \(\R^d\) and \(K^*=4^dK^2\) on \(\T^d\).
\end{proposition}

\begin{proof}
For each real coordinate \(f\) of \(H\) used in
Lemma~\ref{lem:matrix-distance}, we have
\(|f-f_R|\leq\|H-H_R\|\), so \(K_{a,\mathcal R}(f)\leq K\).
On \(\R^d\), apply \eqref{eq:scalar-rect-correction} directly.
On \(\T^d\), apply it to the periodic extension using
\eqref{eq:periodic-moment} and then periodize the correction by
Lemma~\ref{lem:periodic-gj} with \(q=2\). Each scalar
approximant has supremum norm at most \(a^{-1}\log K^*\) and
quadratic oscillation error at most \(C_d/a\).
Assemble the approximants into a Hermitian matrix \(B\) as in
Lemma~\ref{lem:matrix-distance}. Write \(X=(H-B)-(H-B)_R\).
The Frobenius norm bounds the operator norm and satisfies
\[
 \|X\|_{\mathrm F}^2
 =\sum_i|X_{ii}|^2+
 2\sum_{i<j}\bigl(|\operatorname{Re}X_{ij}|^2
                  +|\operatorname{Im}X_{ij}|^2\bigr).
\]
Averaging gives \(\fint_R\|X\|^2\leq(2N^2-N)(C_d/a)^2\).
The same norm identity applied to \(B\) gives the claimed
supremum bound.
\end{proof}

\begin{corollary}\label{cor:rect-log-distance}
If \(d>1\) and \(W\in\mathbf A_p(\mathcal R)\) on \(\R^d\)
or \(\T^d\), \(1\leq p<\infty\), then
\[
 \mathcal D_{\mathcal R}(W)\leq\sqrt{2N^2-N}\,C_d M_p .
\]
\end{corollary}

\begin{proof}
Theorem~\ref{thm:log-oscillation} gives the hypothesis of
Proposition~\ref{prop:rect-correction} at \(a=1/M_p\) with
\(K\leq 2[W]_{\mathbf A_p(\mathcal R)}^{1/M_p}\).
Use Cauchy--Schwarz in \eqref{eq:rect-quadratic}.
\end{proof}

\subsection{A two-sided estimate under John--Nirenberg}
\label{subsec:two-sided}
Scalar John--Nirenberg gives a lower distance bound in terms of
the critical exponential scale. For
\(H\in\mathrm{BMO}_{\mathrm{op}}(\B)\) set
\begin{equation}\label{eq:eta-definition}
 \eta_{\B,N}(H)=\inf\left\{\lambda>0:
 \sup_{E\in\B}\fint_E
 e^{\|H-H_E\|/\lambda}<\infty\right\},
 \qquad\inf\varnothing=+\infty.
\end{equation}
This threshold measures exponential oscillation, rather than
membership of \(e^{H/\lambda}\) in matrix \(A_2\); the two are
compared in Section~\ref{sec:comparison}.

\begin{proposition}\label{prop:eta-distance}
Let \(\B\) have the scalar John--Nirenberg property with constants
\(J,c\). Then, for every Hermitian
\(H\in\mathrm{BMO}_{\mathrm{op}}(\B)\),
\begin{equation}\label{eq:eta-lower-distance}
 \frac c2\,\eta_{\B,N}(H)
 \leq\dist_{\mathrm{BMO}_{\mathrm{op}}(\B)}
 \bigl(H,L^\infty(\Omega;\mathcal H_N)\bigr).
\end{equation}
If \(\B\) also has the scalar Garnett--Jones property with constant
\(G\), then
\begin{equation}\label{eq:eta-two-sided-distance}
 \frac c2\,\eta_{\B,N}(H)
 \leq\dist_{\mathrm{BMO}_{\mathrm{op}}(\B)}
 \bigl(H,L^\infty(\Omega;\mathcal H_N)\bigr)
 \leq(2N^2-N)G\,\eta_{\B,N}(H).
\end{equation}
In particular, the lower-bound constant is independent of \(N\).
\end{proposition}

\begin{proof}
\emph{Step 1: bounded perturbations preserve \(\eta_{\B,N}\).}
If \(F=H-B\) with \(B\) bounded Hermitian, then
\[
 \|(H-H_E)-(F-F_E)\|\leq2\|B\|_\infty,
\]
so the exponential averages for \(H\) and \(F\) at any fixed
\(\lambda>0\) differ by at most the factor \(e^{2\|B\|_\infty/\lambda}\).
Hence \(\eta_{\B,N}(H)=\eta_{\B,N}(F)\).

\emph{Step 2: a scalar majorant.} Put
\(r=\|F\|_{\mathrm{BMO}_{\mathrm{op}}(\B)}\). If \(r=0\), then
\(F=F_E\) a.e.\ on every \(E\in\B\), so \(\eta_{\B,N}(F)=0\).
Suppose \(r>0\). For fixed \(E\in\B\), define the real function
\(u_E(x)=\|F(x)-F_E\|\). For any \(A\in\B\), the reverse triangle
inequality gives
\[
 \fint_A|u_E-(u_E)_A|
 \leq2\fint_A\bigl|u_E-\|F_A-F_E\|\bigr|
 \leq2\fint_A\|F-F_A\|\leq2r.
\]
Thus \(\|u_E\|_{\mathrm{BMO}(\B)}\leq2r\), while
\((u_E)_E\leq r\).

\emph{Step 3: exponential bound.} By Lemma~\ref{lem:layer-cake}
with \(s=2r\), for \(0<a<c/(2r)\),
\[
 \fint_E e^{a\|F-F_E\|}
 \leq e^{ar}\fint_E e^{a|u_E-(u_E)_E|}
 \leq e^{ar}\left(1+\frac{Ja}{c/(2r)-a}\right).
\]
The bound is uniform in \(E\), so
\(\eta_{\B,N}(H)=\eta_{\B,N}(F)\leq2r/c\). Taking the infimum over
\(B\) proves \eqref{eq:eta-lower-distance}.

For every admissible \(\lambda\) in \eqref{eq:eta-definition},
Lemma~\ref{lem:matrix-distance} applies with \(a=1/\lambda\).
Taking the infimum over \(\lambda\) proves the upper bound in
\eqref{eq:eta-two-sided-distance}.
\end{proof}

On a basis with both scalar properties, the John--Nirenberg
estimate also gives an upper bound linear in \(N\). Indeed,
integrating \eqref{eq:scalar-john-nirenberg} bounds the quadratic
oscillation of any real residual by \(\sqrt{2J}/c\) times its
BMO seminorm. Applying this to the coordinate residuals in
Lemma~\ref{lem:matrix-distance} and summing their squares as in
Proposition~\ref{prop:rect-correction} gives the upper constant
\(\sqrt{2N^2-N}\,(\sqrt{2J}/c)G\) in
\eqref{eq:eta-two-sided-distance}. In particular, a bound of this
order holds on cubes. For rectangles the direct correction
estimate gives the constant in the following corollary and
controls the size of the approximant.

On Euclidean rectangles, the scalar exponential estimate in
\cite[Section 9]{NieSik26} gives the scalar John--Nirenberg
property by Markov's inequality. On torus rectangles,
\eqref{eq:periodic-bmo} gives the property by restriction from
the Euclidean estimate. We use common constants \(J_d=B_0^d\)
and \(c_d=2^{-(d+1)}(\log 2)/8\) for the scalar John--Nirenberg
property on both bases, where \(B_0\) is the interval moment
bound in \cite[Section 9]{NieSik26}.

\begin{corollary}\label{cor:rect-eta-distance}
Let \(d>1\) and let \(\mathcal R\) be the rectangle basis on
\(\Omega=\R^d\) or \(\T^d\). For every Hermitian
\(H\in\mathrm{BMO}_{\mathrm{op}}(\mathcal R)\),
\[
 \frac{c_d}{2}\eta_{\mathcal R,N}(H)
 \leq\dist_{\mathrm{BMO}_{\mathrm{op}}(\mathcal R)}
       \bigl(H,L^\infty(\Omega;\mathcal H_N)\bigr)
 \leq\sqrt{2N^2-N}\,C_d\,\eta_{\mathcal R,N}(H).
\]
\end{corollary}

\begin{proof}
Use Proposition~\ref{prop:eta-distance} for the lower bound.
For the upper bound apply Proposition~\ref{prop:rect-correction}
at \(a=1/\lambda\) for each admissible \(\lambda\) in
\eqref{eq:eta-definition}, and take the infimum.
\end{proof}

For \(H=\log W\), Theorem~\ref{thm:log-oscillation} gives
\(\eta_{\B,N}(\log W)\leq M_p\). The two-sided distance
estimates apply to such logarithms on cubes by
Proposition~\ref{prop:eta-distance} and on rectangles by
Corollary~\ref{cor:rect-eta-distance}.

\section{Necessity of the scalar Garnett--Jones property}
\label{sec:equivalence}

We prove that finiteness of \(L_{\B,N}\) is equivalent to the
scalar Garnett--Jones property. Under the small-BMO condition,
both are equivalent to scalar \(A_2\) decomposition.
Recall that \(L_{\B,1}=\sup_{w\in A_2(\B)}
\dist_{\mathrm{BMO}(\B)}(\log w,L^\infty)\).

The following elementary lemma allows us to test matrix distance
bounds on scalar weights: embedding \(w\) as \(wI_N\) preserves
both the \(A_2\) characteristic and the distance to bounded functions.

\begin{lemma}\label{lem:scalar-embedding}
For every \(w\in A_2(\B)\) and \(N\geq1\),
\([wI_N]_{\mathbf A_2(\B)}=[w]_{A_2(\B)}\) and
\begin{equation}\label{eq:scalar-matrix-distance}
 \mathcal D_{\B}(wI_N)
 =\dist_{\mathrm{BMO}(\B)}(\log w,L^\infty).
\end{equation}
\end{lemma}

\begin{proof}
The identity of characteristics is immediate from
\eqref{eq:ap-definition}. Approximants \(bI_N\) with \(b\) bounded
real have operator-norm oscillation equal to the scalar oscillation
of \(\log w-b\), which gives ``\(\leq\)'' in
\eqref{eq:scalar-matrix-distance}. Conversely, for any bounded
Hermitian \(B\), the entry \(B_{11}\) is bounded and real, and
\eqref{eq:entrywise} gives
\(\|\log w-B_{11}\|_{\mathrm{BMO}(\B)}
\leq\|(\log w)I_N-B\|_{\mathrm{BMO}_{\mathrm{op}}(\B)}\).
\end{proof}

\begin{theorem}\label{thm:gj-equivalence}
Let \(\B\) be an averaging basis and \(G>0\). The following are
equivalent:
\begin{enumerate}
\item \(\B\) has the scalar Garnett--Jones property with constant
\(G\);
\item every \(w\in A_2(\B)\) satisfies
\(\dist_{\mathrm{BMO}(\B)}(\log w,L^\infty)\leq G\).
\end{enumerate}
Moreover, as extended real numbers,
\begin{equation}\label{eq:optimal-distance-constants}
 G_{\B}^{\mathrm{opt}}=L_{\B,1}
 \leq L_{\B,N}
 \leq(2N^2-N)G_{\B}^{\mathrm{opt}}
 \qquad(N\geq1).
\end{equation}
In particular, \(L_{\B,N}<\infty\) for some, equivalently all,
\(N\geq1\) if and only if \(\B\) has the scalar Garnett--Jones
property.
\end{theorem}

\begin{proof}
\emph{(1)\(\Rightarrow\)(2).} Apply Theorem~\ref{thm:log-distance}
with \(N=1\) and \(p=2\), where \(2N^2-N=1\) and \(M_2=1\).

\emph{(2)\(\Rightarrow\)(1).} Let \(f\in\mathrm{BMO}(\B)\) be real,
\(a>0\), and \(K=\sup_{E\in\B}\fint_Ee^{a|f-f_E|}<\infty\). Put
\(w=e^{af}\). For \(E\in\B\), \(f_E\) is finite and
\[
 \fint_Ew=e^{af_E}\fint_Ee^{a(f-f_E)}\leq e^{af_E}K,
 \qquad
 \fint_Ew^{-1}\leq e^{-af_E}K.
\]
So \(w^{\pm1}\) are integrable on basis sets, hence locally
integrable by finite coverings, and \([w]_{A_2(\B)}\leq K^2\). By
(2) and \(\log w=af\),
\[
 \dist_{\mathrm{BMO}(\B)}(f,L^\infty)
 =\frac1a\dist_{\mathrm{BMO}(\B)}(\log w,L^\infty)
 \leq\frac Ga.
\]
The independence of the bound in (2) from \([w]_{A_2(\B)}\) is
essential here, since \(K^2\) may be arbitrarily large. This
direction uses the assumed exponential-integrability bound
directly and needs no John--Nirenberg inequality for \(\B\).

\emph{The constants.} By the equivalence just proved, \(G>0\) is
admissible in (1) if and only if \(G\geq L_{\B,1}\). Hence
\(G_{\B}^{\mathrm{opt}}=L_{\B,1}\), also when the common value is
\(0\) or \(+\infty\). Lemma~\ref{lem:scalar-embedding} shows that
\(wI_N\) is admissible in \(L_{\B,N}\) with the same distance, so
\(L_{\B,1}\leq L_{\B,N}\). Finally, Theorem~\ref{thm:log-distance}
at \(p=2\) gives \(L_{\B,N}\leq(2N^2-N)G\) for every admissible
\(G\), and taking the infimum gives the last inequality.
\end{proof}

\subsection{\texorpdfstring{The \(A_2\) decomposition property}{The A2 decomposition property}}
\label{subsec:decomposition}

\begin{definition}[{\cite[Definition 3.1]{NieSik19}}]
\label{def:a2-decomposition}
The basis \(\B\) has the \emph{\(A_2\)-decomposition property} if
there are \(0<\delta\leq1\) and \(K_{\mathrm d}\geq1\) such that
every \(w\in A_2(\B)\) can be written
\begin{equation}\label{eq:a2-decomposition}
 w^\delta=u\varphi,\qquad
 [u]_{A_2(\B)}\leq K_{\mathrm d},\qquad
 \varphi>0,\quad\varphi,\varphi^{-1}\in L^\infty(\Omega).
\end{equation}
Here \(\delta\) and \(K_{\mathrm d}\) are independent of \(w\);
the bounds on \(\varphi^{\pm1}\) may depend on \(w\).
\end{definition}

To pass from a distance estimate back to a decomposition, we use a
uniform small-BMO condition: there are \(\beta>0\) and
\(K_0\geq1\) such that, for every real \(h\in\mathrm{BMO}(\B)\),
\begin{equation}\label{eq:small-bmo-exponential}
 \|h\|_{\mathrm{BMO}(\B)}\leq\beta
 \quad\Longrightarrow\quad
 [e^h]_{A_2(\B)}\leq K_0.
\end{equation}
This is an additional assumption on \(\B\); by
Remark~\ref{rem:small-bmo-jn} it follows from the scalar
John--Nirenberg property.

\begin{theorem}\label{thm:decomposition-equivalence}
Let \(\B\) be an averaging basis.
\begin{enumerate}
\item If \(\B\) has the \(A_2\)-decomposition property with
parameters \(\delta,K_{\mathrm d}\), then it has the scalar
Garnett--Jones property with constant
\begin{equation}\label{eq:decomposition-to-gj}
 G=\frac{\log(2K_{\mathrm d})}{\delta}.
\end{equation}
\item If \(\B\) satisfies \eqref{eq:small-bmo-exponential} and has
the scalar Garnett--Jones property with constant \(G\), then it
has the \(A_2\)-decomposition property with
\begin{equation}\label{eq:gj-to-decomposition}
 \delta=\min\Bigl\{1,\frac{\beta}{2G}\Bigr\},\qquad
 K_{\mathrm d}=K_0.
\end{equation}
\end{enumerate}
Consequently, under \eqref{eq:small-bmo-exponential}, the
\(A_2\)-decomposition property, the scalar Garnett--Jones property,
and finiteness of \(L_{\B,N}\) for any fixed \(N\geq1\) are
equivalent.
\end{theorem}

\begin{proof}
(1) Let \(w\in A_2(\B)\) and write \(w^\delta=u\varphi\) as in
\eqref{eq:a2-decomposition}. Then \(\delta^{-1}\log\varphi\) is
bounded and real, and
\(\log w-\delta^{-1}\log\varphi=\delta^{-1}\log u\).
By \eqref{eq:log-bmo} for the scalar weight \(u\) at \(p=2\),
\[
 \dist_{\mathrm{BMO}(\B)}(\log w,L^\infty)
 \leq\frac{\|\log u\|_{\mathrm{BMO}(\B)}}{\delta}
 \leq\frac{\log[u]_{A_2(\B)}+\log2}{\delta}
 \leq\frac{\log(2K_{\mathrm d})}{\delta}.
\]
Theorem~\ref{thm:gj-equivalence} gives the scalar Garnett--Jones
property with this constant.

(2) Let \(w\in A_2(\B)\). By Theorem~\ref{thm:gj-equivalence},
\(\dist(\log w,L^\infty)\leq G<2G\), so there is a bounded real
\(g\) with \(\|\log w-g\|_{\mathrm{BMO}(\B)}\leq2G\). With
\(\delta\) as in \eqref{eq:gj-to-decomposition}, put
\[
 h=\delta(\log w-g),\qquad u=e^h,\qquad\varphi=e^{\delta g}.
\]
Then \(\|h\|_{\mathrm{BMO}(\B)}\leq2\delta G\leq\beta\), so
\([u]_{A_2(\B)}\leq K_0\) by \eqref{eq:small-bmo-exponential};
\(\varphi^{\pm1}\) are bounded; and \(u\varphi=w^\delta\).

The final assertion follows from (1), (2), and
Theorem~\ref{thm:gj-equivalence}.
\end{proof}

Theorems~\ref{thm:gj-equivalence} and
\ref{thm:decomposition-equivalence} give Theorem~\ref{thmC}.

\begin{remark}\label{rem:small-bmo-jn}
If \(\B\) has the scalar John--Nirenberg property with constants
\(J,c\), then \eqref{eq:small-bmo-exponential} holds with
\(\beta=c/2\) and \(K_0=(1+J)^2\). Indeed, if
\(\|h\|_{\mathrm{BMO}(\B)}\leq c/2\), Lemma~\ref{lem:layer-cake}
with \(s=c/2\) and \(a=1\) gives
\(\fint_Ee^{|h-h_E|}\leq1+J\), and therefore
\[
 \Bigl(\fint_Ee^h\Bigr)\Bigl(\fint_Ee^{-h}\Bigr)
 =\Bigl(\fint_Ee^{h-h_E}\Bigr)\Bigl(\fint_Ee^{h_E-h}\Bigr)
 \leq(1+J)^2;
\]
finite coverings give local integrability of \(e^{\pm h}\). Thus
the equivalences in Theorem~\ref{thm:decomposition-equivalence}
hold on every basis with the scalar John--Nirenberg property.
\end{remark}

The scalar rectangle theorem settles the decomposition problem
for Euclidean rectangles left open in \cite{NieSik19}.
The torus transfer gives the same property on \(\T^d\).

\begin{corollary}\label{cor:rect-decomposition}
Let \(d>1\) and let \(\mathcal R\) be the rectangle basis on
\(\R^d\) or \(\T^d\). Every scalar \(w\in A_2(\mathcal R)\) has
the decomposition \eqref{eq:a2-decomposition} with
\[
 \delta=\min\{1,c_d/(4C_d)\},
 \qquad K_{\mathrm d}=(1+J_d)^2,
\]
where \(J_d,c_d\) are the scalar John--Nirenberg constants for
rectangles used above.
\end{corollary}

\begin{proof}
Use \(G=C_d\) from \eqref{eq:scalar-rect-correction} and
Lemma~\ref{lem:periodic-gj}. Remark~\ref{rem:small-bmo-jn}
gives \(\beta=c_d/2\) and \(K_0=(1+J_d)^2\), so the result
follows from Theorem~\ref{thm:decomposition-equivalence}.
\end{proof}

For the comparison with \cite{NieSik19}, define for
\(H\in L^1_{\mathrm{loc}}(\Omega;\mathcal H_N)\) the matrix
\(A_2\) threshold
\begin{equation}\label{eq:epsilon-definition}
 \varepsilon_{\B,N}(H)
 =\inf\{\lambda>0:e^{H/\lambda}\in\mathbf A_2(\B)\},
 \qquad\inf\varnothing=+\infty,
\end{equation}
and write \(\varepsilon_\B=\varepsilon_{\B,1}\) for real \(f\);
this is the threshold of \cite{NieSik19}. The scalar Garnett--Jones
property can be expressed through this threshold without changing
the constant.

\begin{proposition}\label{prop:epsilon-upper-distance}
For \(G>0\), \(\B\) has the scalar Garnett--Jones property with
constant \(G\) if and only if
\begin{equation}\label{eq:epsilon-upper-distance}
 \dist_{\mathrm{BMO}(\B)}(f,L^\infty)
 \leq G\,\varepsilon_{\B}(f)
 \quad\text{for every real }f\in\mathrm{BMO}(\B).
\end{equation}
\end{proposition}

\begin{proof}
Assume the scalar Garnett--Jones property with constant \(G\). If
\(e^{f/\lambda}\in A_2(\B)\), Theorem~\ref{thm:gj-equivalence}
applied to this weight gives
\(\dist(f,L^\infty)=\lambda\dist(f/\lambda,L^\infty)\leq G\lambda\).
Taking the infimum over \(\lambda\) gives
\eqref{eq:epsilon-upper-distance}; if no \(\lambda\) is admissible
the right-hand side is \(+\infty\).

Conversely, assume \eqref{eq:epsilon-upper-distance}, and let
\(a>0\) with \(\sup_E\fint_Ee^{a|f-f_E|}<\infty\). As in the proof
of Theorem~\ref{thm:gj-equivalence}, \(e^{af}\in A_2(\B)\), so
\(\varepsilon_\B(f)\leq1/a\) and \(\dist(f,L^\infty)\leq G/a\).
\end{proof}

\begin{remark}
In \cite[Theorem 3.5]{NieSik19} the equivalence between the
decomposition property and \eqref{eq:epsilon-upper-distance} is
proved on \(\R^d\) for structured Muckenhoupt bases with the
Calder\'on--Zygmund decomposition property, and
\cite[Lemma 3.4]{NieSik19} supplies exactly
\eqref{eq:small-bmo-exponential} for such bases.
Proposition~\ref{prop:epsilon-upper-distance} identifies
\eqref{eq:epsilon-upper-distance} with the scalar Garnett--Jones
property without loss of constants, and
Theorem~\ref{thm:decomposition-equivalence} shows that
\eqref{eq:small-bmo-exponential} is the only additional input
needed for the converse.
\end{remark}

\section{\texorpdfstring{Sharpness of the estimate in \(p\)}{Sharpness of the estimate in p}}\label{sec:sharpness}

Let \(W\in\mathbf A_p(\B)\). Spectral truncations of \(\log W\) are
bounded and converge to \(\log W\) in \(L^q_{\mathrm{loc}}\), but in
general not in \(\mathrm{BMO}_{\mathrm{op}}(\B)\). For \(R>0\) let
\(\tau_R(t)=\max\{-R,\min\{t,R\}\}\) and
\(H_R=\tau_R(\log W)\), defined by functional calculus. Then
\(\|H_R\|_\infty\leq R\) and
\(\|\log W-H_R\|=(\|\log W\|-R)_+\leq\|\log W\|\), so
Theorem~\ref{thm:log-oscillation} and dominated convergence give
\(H_R\to\log W\) a.e.\ and in \(L^q_{\mathrm{loc}}\) for every
\(1\leq q<\infty\). Convergence in \(\mathrm{BMO}_{\mathrm{op}}\)
would force \(\mathcal D_\B(W)=0\). The following power weights
show that the distance is positive, and that its size can be of
order \(M_p\).

\begin{proposition}\label{prop:positive-distance}
On \(\R\), let \(W(x)=|x|^\alpha I_N\), with \(W(0)\) an arbitrary
positive definite matrix. If \(1<p<\infty\), assume
\(-1<\alpha<p-1\); if \(p=1\), assume \(-1<\alpha\leq0\). Then
\(W\in\mathbf A_p(\mathcal Q)\) and
\[
 \mathcal D_{\mathcal Q}(W)\geq\frac{|\alpha|\log2}{2}.
\]
\end{proposition}

\begin{proof}
For \(W=wI_N\) one has
\(\|W^{1/p}(x)W^{-1/p}(y)\|=(w(x)/w(y))^{1/p}\), so the matrix
characteristics equal the scalar ones. For \(w(x)=|x|^\alpha\), the
stated ranges of \(\alpha\) are exactly those for which
\(w\in A_p(\R)\), respectively \(w\in A_1(\R)\); see
\cite[Example 7.1.7]{Gra14}.

For the lower bound, let \(B\) be bounded Hermitian, \(b=B_{11}\),
\(f(x)=\alpha\log|x|\), \(g=f-b\), and \(I_k=(0,2^{-k})\). Then
\(f_{I_k}=-\alpha(k\log2+1)\), so
\(|g_{I_k}-g_{I_0}|\geq k|\alpha|\log2-2\|b\|_\infty\). Since
\(|I_j|=2|I_{j+1}|\),
\[
 |g_{I_{j+1}}-g_{I_j}|
 \leq\fint_{I_{j+1}}|g-g_{I_j}|
 \leq2\|g\|_{\mathrm{BMO}(\R)}.
\]
Summing over \(j=0,\ldots,k-1\), dividing by \(k\), and letting
\(k\to\infty\) gives \(\|g\|_{\mathrm{BMO}}\geq|\alpha|\log2/2\).
By \eqref{eq:entrywise},
\(\|\log W-B\|_{\mathrm{BMO}_{\mathrm{op}}}\geq\|g\|_{\mathrm{BMO}}\).
Take the infimum over \(B\).
\end{proof}

The supremum of \(|\alpha|\) over the admissible range is
\(\max\{1,p-1\}=M_p\) for every \(1\leq p<\infty\). Combined with
Corollary~\ref{cor:cube-log-distance}, this gives the following.

\begin{corollary}\label{cor:sharp-p}
Let \(\mathcal Q\) be the interval basis on \(\R\), let
\(1\leq p<\infty\) and \(N\geq1\), and put
\(L^{(p)}_{\mathcal Q,N}=\sup\{\mathcal D_{\mathcal Q}(W):
W\in\mathbf A_p(\mathcal Q),\ W\text{ is }N\times N\}\). Then
\[
 \frac{\log2}{2}\,M_p
 \leq L^{(p)}_{\mathcal Q,N}
 \leq(2N^2-N)\,G_1\,M_p .
\]
In particular, the supremum grows linearly in \(p\) as
\(p\to\infty\) and has a uniform positive lower bound for
\(1\leq p\leq2\).
\end{corollary}

For Euclidean rectangles in dimension \(d>1\), the same lower
bound follows by taking \(W(x)=|x_1|^\alpha I_N\).
Averaging on a product shows that the characteristic equals
the interval characteristic of \(|x_1|^\alpha\). For a bounded
Hermitian approximant, use the first diagonal entry and the
rectangles \(I_k\times(0,1)^{d-1}\) in the proof of
Proposition~\ref{prop:positive-distance}. Their successive
measure ratios are two, and the averages of
\(\alpha\log|x_1|\) are the same as on \(I_k\). It follows that
\[
 \frac{\log2}{2}M_p
 \leq
 \sup_{W\in\mathbf A_p(\mathcal R),\ W\text{ is }N\times N}
      \mathcal D_{\mathcal R}(W)
 \leq \sqrt{2N^2-N}\,C_d M_p .
\]

\section{\texorpdfstring{The matrix \(A_2\) threshold and comparison with earlier results}{The matrix A2 threshold and comparison with earlier results}}
\label{sec:comparison}

The comparison in \cite[Proposition 4.5]{NieSik21} is for
cubes on \(\T^d\). Its extension in
\cite[Proposition 17.1]{NieSik25} is for structured torus
bases satisfying the geometric condition (17.9) of that work.
The constants there can depend on the prescribed log-preserving
unitary matrix function \(U(x)\). Conjugation by a variable
\(U(x)\) need not commute with averaging. We use
Theorem~\ref{thm:log-distance} to bound the distance taken over
all bounded Hermitian approximants in operator-norm BMO.

\begin{corollary}\label{cor:matrix-epsilon}
Let \(\B\) have the scalar Garnett--Jones property with constant
\(G\), and let \(H\in L^1_{\mathrm{loc}}(\Omega;\mathcal H_N)\).
Then
\begin{equation}\label{eq:matrix-epsilon-distance}
 \dist_{\mathrm{BMO}_{\mathrm{op}}(\B)}
       \bigl(H,L^\infty(\Omega;\mathcal H_N)\bigr)
 \leq(2N^2-N)\,G\,\varepsilon_{\B,N}(H).
\end{equation}
\end{corollary}

\begin{proof}
If \(\varepsilon_{\B,N}(H)=+\infty\) there is nothing to prove.
For admissible \(\lambda\), Theorem~\ref{thm:log-distance} at
\(p=2\) applied to \(e^{H/\lambda}\) bounds the distance of
\(H/\lambda\) by \((2N^2-N)G\), and the distance is positively
homogeneous. Take the infimum over \(\lambda\).
\end{proof}

The matrix \(A_2\) threshold \(\varepsilon_{\B,N}(H)\) must be
distinguished from the exponential-oscillation threshold
\(\eta_{\B,N}(H)\) of \eqref{eq:eta-definition}. For
\(H\in\mathrm{BMO}_{\mathrm{op}}(\B)\),
Theorem~\ref{thm:log-oscillation} applied to \(e^{H/\lambda}\) at
\(p=2\) gives
\[
 \eta_{\B,N}(H)\leq\varepsilon_{\B,N}(H).
\]
For \(N=1\), equality holds: finite exponential oscillation at
scale \(\lambda\) implies \(e^{H/\lambda}\in A_2(\B)\) by bounding
the averages of \(e^{\pm(H-H_E)/\lambda}\). In the matrix case,
\eqref{eq:eta-two-sided-distance} does not by itself yield a lower
distance bound in terms of \(\varepsilon_{\B,N}(H)\), since the
comparison of thresholds goes in the opposite direction.
For Euclidean and torus rectangles in dimension \(d>1\),
Corollary~\ref{cor:rect-eta-distance} and the threshold comparison
give \eqref{eq:matrix-epsilon-distance} with
\(\sqrt{2N^2-N}\,C_d\) as the constant.

On the circle with its interval basis \(\mathcal Q\), Bloom's
example \cite[pp.~100--103]{Blo81} gives a Hermitian BMO function
\(H\) with \(e^{tH}\notin\mathbf A_2\) for every \(t>0\).
This holds also for \(t<0\), since taking adjoints and exchanging
\(x,y\) in \eqref{eq:ap-definition} at \(p=2\) shows that
\(W\mapsto W^{-1}\) preserves \(\mathbf A_2\) and its characteristic.
The scalar John--Nirenberg property and
Proposition~\ref{prop:eta-distance} give
\(\eta_{\mathcal Q,N}(H)<\infty\), whereas
\(\varepsilon_{\mathcal Q,N}(H)=+\infty\), so
\eqref{eq:matrix-epsilon-distance} is vacuous. Small matrix BMO
does not in general imply matrix \(A_2\) exponentiation.

\section{Open questions}\label{sec:questions}

The lower bounds for \(L_{\B,N}\) from scalar weights are
independent of \(N\). On cubes and rectangles the upper bounds
in Section~\ref{sec:distance} are linear in \(N\). Can
\(L_{\B,N}\) be bounded independently of \(N\) on rectangles
or on cubes? It would also be of interest to find a basis
satisfying the scalar John--Nirenberg property for which the
scalar Garnett--Jones estimate fails. By
Theorem~\ref{thm:decomposition-equivalence} and
Remark~\ref{rem:small-bmo-jn}, this is equivalent to finding a
John--Nirenberg basis where scalar \(A_2\) decomposition fails.
Finally, in view of Bloom's example, a lower bound for the matrix
distance in terms of \(\varepsilon_{\B,N}\), rather than
\(\eta_{\B,N}\), would require an additional structural
assumption on \(H\).

\section*{Acknowledgment}
ChatGPT Sol 6.1 was used under the authors' direction in developing the
construction. AI tools also assisted with proof review and editing.

\end{document}